\input amstex
\documentstyle{amsppt}
\NoBlackBoxes
\topmatter
\title Self Coincidence Numbers and the Fundamental Group
\endtitle
\author Daniel Henry Gottlieb 
\endauthor
\document
\bigskip
\dedicatory I am very happy to dedicate this paper to my two friends, Ed Fadell and Albrecht Dold on the occasion of their 80th birthdays, and I thank them for being older than I since I was able to use their
powerful theorems as a young mathematician.--- Dan Gottlieb
\enddedicatory
\abstract For $M$ and $N$ closed oriented connected smooth manifolds of the same dimension, we consider the mapping space ${\text Map}(M,N;f)$ of continuous maps homotopic to $f:M\rightarrow N$.
We will show that the evaluation map from the space of maps to the manifold $N$ induces a nontrivial homomorphism on the fundamental group only if the self coincidence number of $f$, denoted 
$\Lambda_{f,f}$, equals zero. Since $\Lambda_{f,f}$ is equal to the product of the degree of $f$ and the Euler--Poincare number of $N$, we obtain  results related to earlier results about the evaluation map and the Euler--Poincare number. \endabstract 
\endtopmatter
\medskip
\centerline{\bf 1. INTRODUCTION}

\bigskip\bigskip

This paper owes its existence to a question and conjecture of Dusa McDuff. McDuff was preparing her
paper [McDuff (2006)] on the Symplectomorphism group of a blow up. She wanted to prove a variation
of my theorem that the evaluation map $\omega : {\text{Map}(X,X;1_X)} \rightarrow X$ induces a trivial homomorphism on the fundamental group if $X$ is compact CW complex whose Euler--Poincare number
is not zero; [Gottlieb(1965), Theorem II.7]. McDuff's conjecture was that $\omega : {\text{Map}(M,N;f)} \rightarrow N :  g \mapsto \omega(g) := g(*)$ induced the trivial homomorphism on the fundamental 
groups in the case where $N$ is a smooth orientable manifold and $M$ is a blow up.

McDuff asked me if there was a different proof of my original result besides the original proof using
Neilsen-Wecken theory. There was one, in [Gottlieb (1990); see section 9]. I told her that and included
suggestions for proceeding using coincidence theory. After awhile I thought about whether my suggested method would actually work, and I began thinking about it in more detail, realizing that I needed to use self-coincidence numbers. This resulted in Theorem 1.1 below. However, McDuff managed to come up with an independent proof for the blow up case in [McDuff (2006); see 
Lemma 2.2 iii].

\proclaim{Theorem 1.1} Let $f: M^n \rightarrow N^n$ be a map of closed oriented connected smooth manifolds of the same dimension. Let $\omega : {\text{Map}(M,N;f)} \rightarrow N :  g \mapsto \omega(g) 
:= g(*)$
where $*$ is a base point of $M$. Then if $\ \omega_* : \pi_1(\text{Map}(M,N;f) \rightarrow \pi_1(N)$ is non-trivial, either the degree of $f$, $deg(f) =0$ or the Euler Characteristic $\chi (N)=0$.
\endproclaim

Note Theorem 1.1 does imply McDuff's Lemma 2.2(iii) since a blow up $f$ has degree $f$ equal to one.

Now assuming that $N$ is aspherical implies that the fundamental group of the mapping space can be
expressed in group theoretic terms and we have the following corollaries, which are related to the 
result in [Gottlieb (1965)] which states that the center of the fundamental group of  compact aspherical 
complex is trivial if the Euler--Poincare number is not zero.

\proclaim{Corollary 1.2} In addition to the hypotheses of Theorem 1.1 for $M$ and $N$, we assume that $N$ is aspherical.  Then if $deg(f)\times \chi(N) \neq 0$ , then  $\text{Map}(M,N;f)$ is contractible.
\endproclaim

\proclaim{Corollary 1.3} Let $N$ be a closed oriented connected aspherical manifold with 
$ \chi(N) \neq 0$. Then every subgroup of finite index in $\pi_1(N)$ has a trivial centralizer.
\endproclaim

Finally, Theorem 1.1 follows from theorem 1.4 and a calculation of the self coincidence number.

\proclaim{Theorem 1.4} Let $f: M^n \rightarrow N^n$ be a map of closed oriented connected smooth manifolds of the same dimension. Let $\omega : \text{Map}(M,N;f) \rightarrow N :  g \mapsto \omega(g) 
:= g(*)$
where $*$ is a base point of $M$. Then if  $\ \omega_* : \pi_1(\text{Map}(M,N;f) \rightarrow \pi_1(N)$ is non-trivial, the Lefchetz self coincidence number $ \Lambda_{f,f} = 0$.
\endproclaim

There is considerable interest in the self--coincidence theory at this time, see [Dold and Goncalves (2005)] and [Ulrich Koschorke (2006)].

The above theorems will be proven in Section 3. In Section 2, the key lemmas will be stated and proved
using  new axioms for the coincidence index, due to Christopher Staecker. In Section 4, we will propose
a conjecture vastly generalizing Corollary 1.2, and prove it is true in the case that $N$ is $S^1$.


\bigskip\bigskip
\centerline{\bf 2. THE COINCIDENCE INDEX}

\bigskip\bigskip

In this section, we use the the axioms of Chris Staecker which characterize the coincidence index on a differential manifold. We summarize his results here. See [Staecker (2006)].

For two mappings $f,g : X \rightarrow Y $ we say that $ x \in X$ is a coincidence point of $f$ and $g$ if 
$f(x)=g(x)$. We write $x \in \roman {Coin}(f,g)$

For continuous maps $f,g : X \rightarrow Y$ and some open set $U \subset X$, we say that the triple
$(f,g,U)$ is { \it admissible} if the set of coincidences of $f$ and $g$ in $U$ is compact. Let $\Cal C$ be the set of admissible triples of mappings where $X$ and $Y$ are oriented differentiable manifolds of the same dimension.

If $f_t, g_t: X\times I \rightarrow Y$ are homotopies and $U \subset X$ is an open subset, we say that 
$(f_t,g_t)$ is a pair of {\it admissible homotopies in} $U$ if 
$$\{(x,t)\mid x \in \roman {Coin}(f_t,g_t) \cap U\}$$ is a compact subset of $X\times I$. In this case we say that $(f_0,g_0, U)$ is {\it admissibly homotopic} to $(f_1,g_1, U)$.
The coincidence index is the unique function $\roman{ind}:\Cal C \rightarrow \Bbb R$ which satisfies the following three axioms:

\proclaim {\bf Axiom 1 (Normalization)} Let $f, g:X\rightarrow Y$ be smooth mappings, and  $ \Lambda_{f,g}$ be the Lefschetz coincidence number for $f$ and $g$. Then 
$\roman{ind}(f, g, X)=  \Lambda_{f,g} $
\endproclaim

\proclaim {\bf Axiom 2 (Additivity)} Given maps $f,g:X\rightarrow Y$ and an admissible triple $(f,g,U)$ , if $U_1$ and $U_2$ are disjoint open subsets of $U$ with $\roman{Coin}(f,g) \subset U_1 \cup U_2$, then
$$ \roman {ind}(f,g, U)= \roman{ind}(f,g, U_1)+\roman{ind}(f,g, U_2).$$
\endproclaim

\proclaim{\bf Axiom 3 (Homotopy)} If there exists an admissable homotopy of the triple $(f_0,g_0,U)$ to
$(f_1,g_1,U)$ , then
$$\roman{ind}(f_0,g_0,U)= \roman{ind}(f_1,g_1,U).$$

\endproclaim

\proclaim{\bf Lemma 2.1}
Suppose we have a homotopy $H:X\times I \rightarrow Y$ so that there is a {\it connected} open set
$V \subset X \times I $ which contains a compact subset of $ \roman{Coin}(f_t,g_t) \cap V$. Then if
$V_0:= V\cap X \times 0$ is not empty and if $(f,g,V_0)$ is an admissible triple with $\roman {ind}(f,g, V_0)=: k \neq 0$, then $ V_1:= V\cap X \times 1$ is not empty and  $\roman {ind}(f_1,g_1, V_1)=: k \neq 0$.

\endproclaim

\demo{proof}
Let $V_t := V \cap X\times t$. Assume that $T$ is the first $t$ such that $\roman {ind}(f_t,g_t, V_t) \neq k$. Then about every connected component of  $\roman{Coin}(f_t,g_t) \cap V_T$ we can find open neighborhoods $U$ such that $U \subset V_T$ so that the tube $U\times J \subset V$ where $J$ is a sufficiently small interval about $T$ so that the homotopy restricts on the tube to an admissible homotopy. Since the admissible homotopy on the tubes begins with
$\roman {ind}(f_s,g_s, U)= k$ at $s<T$ , and
it ends where $\roman {ind}(f,g, U)=k$ at a $ t > T$, then $T$ cannot be the first time where 
$\roman {ind}(f,g, V_t)\neq k$. Hence  $\roman {ind}(f_1,g_1, V_1)=: k \neq 0$.
\qed
\enddemo

\proclaim{\bf Lemma 2.2}
Suppose we have homotopies $F:X\times I \rightarrow Y$ and $G:X\times I \rightarrow Y$ associated with $f_t$ and $g_t$. Suppose they are cyclic homotopies so that $f_1 = f_0$ and $g_1=g_0$. Then we can associate them with maps $F:X\times S^1 \rightarrow Y$ and $G : X \times S^1\rightarrow Y$. Assume in addition that there is a closed loop $\sigma$ in $X\times S^1$ such that $f_t(\sigma(t))=g_t(\sigma(t))$. Then 

$$ \omega^f_*(\hat F) + f_*(\sigma_x) = \omega^g_*(\hat G) + g_*(\sigma_x) \in \pi_1(Y) \tag 2.1$$

 where the notation of the conclusion will be explained in the proof.
\endproclaim
 \demo{proof} The cyclic homotopies $F:X\times S^1 \rightarrow Y$ and $G:X\times S^1\rightarrow Y$ are associated  to maps $\hat {F}:S^1 \to \text{Map}(X,Y;f)$ and to $\hat {G}:S^1 \to \text{Map}(X,Y;f)$ respectively. Then we have the evaluation maps $\omega^f: \text{Map}(X,Y;f) \to Y$ and $\omega^g: \text{Map}(X,Y;g) \to Y$. Now $\pi_1 (X\times S^1) \cong \pi_1(X)\times \pi_1(S^1)$ and any element in that fundamental group can be uniquely written in the form $x+k\iota$ where $k$ is an integer and $\iota$ is the generator of
 $\pi_1(S^1)$ and the ``$+$" sign denotes multiplication in the fundamental group, and it reminds us
 that any $x \in \pi_1(X)$ commutes with $\iota$. Now let us use the same symbol for a loop and its 
 homotopy class in the fundamental group. Then $\sigma=: \sigma_x + \iota$ where $\sigma_x$ is the element of $\pi_1(X)$ defined by $\sigma$.
 Now  $F_*(\sigma)=G_*(\sigma) \in \pi_1(Y)$ by the hypothesis on $\sigma$. Also
  $F_*(\sigma)= f_*(\sigma_x) + \omega^f_*(\hat F)$ where we note $F_*(\iota)=\omega^f_*(\hat F)$;
  and also $G$ gives a similar formula. Combining these two formulas yields  the conclusion of the lemma. \qed
  \enddemo
  
  We will be interested in the special case of equation (2.1) when the homotopy $g_t =g$ for all $t$. Thus
  $\omega^g_*(\hat G) = 0$ and we obtain
  $$ \omega^f_*(\hat F) + f_*(\sigma_x) =  g_*(\sigma_x) \in \pi_1(Y) \tag 2.2$$
    
We will be concerned with showing that such a closed loop exists as in the hypotheses of Lemma 2.2.
This requires using Lemma 2.1 and a few extra hypotheses. From Lemma 2.1, if there is a connected open set $V$ containing a coincidence set of points $C:=\{(x,t) \in X \times I : f_t(x)= g(x)\}$ and if that set $C$ 
is path connected, then there is a path $c(s)$ so that $c(0) \in C$ and $f_0(c(0)) =g(c(0))$ and also
$c(1) \in C$ and $f_1(c(1)) = g(c(1))$. Now it is possible that two or more paths of coincidence points
 whose indices add up to zero will annihilate themselves at some time $T$ during the homotopy $f_t$
 so that the path $c(s)$ cannot move forward in homotopy time $t$. But we can reverse the homotopy
 so that  $f_t(c(t)) = g(c(t))$ remains true. Thus by reparameterizing the homotopy we will have
 $f_t(c(t)) = g(c(t))$.
 
 Now one slight problem is that the set $C:=\{(x,t) \in X \times I: f_t(x)= g(x)\}$ may not be a path connected subset.
But in that case we can find a path arbitrarily close to $C$, and thus this path will be close to the set $C$ of coincidences and will satisfy the homotopy condition in Lemma 2.2 so that  equation (2.2) will  hold.

A more substantial difficulty is that the path $c(t)$ of coincidences may not be a closed loop, i.e. 
$c(0)\neq c(1)$ even when the homotopy is cyclic, $f_0 =f_1$. An example of this occurs in the case
that $g: X \to Y$ is a finite covering map and $f_t$ is a constant map for all $t$ and the homotopy traces
out a loop in $Y$ whose lift is not a closed loop in $X$.

One idea of avoiding this problem is to consider the mappings 
$$ X \times S^1\rightarrow X\times Y: (x,t)\mapsto (x,f_t(x))$$
and 
$$X \times S^1\to X\times Y: (x,t) \mapsto  (x,g(x)).$$ Isotopying the second map so it is transverse to the first map
gives rise to an inverse image of the coincidence set in $Y$ as a closed one dimensional manifold in 
$X\times S^1$ of points which are close to the coincidence points $C$. At least one path componant
of this inverse image will be an embedded circle in $X\times S^1$.This will give us a closed loop 
$\sigma$, but $\sigma = \sigma_x + k\iota$ for some integer $k$, not necessarily equal to one.

One way to guarantee closed loops $\sigma_x$ and $\sigma$ is to assume that $\text{Coin}(f, g)$ 
consists of one point, whose coincidence index is not zero. This implies that there is a connected subset of
$C:=\{(x,t) \in X\times I: f_t(x)= g(x)\}$ which connects with the unique coincidence point at $t=1$. Thus
the path $c(t)$ we have constructed will start at the unique coincidence point and pass through $X\times S^1$ and close at $t=1$ to the unique coincidence point.

However, the requirement of only one coincidence point is a strong assumption, which may not be satisfied. This assumption is  major subject of study in Nielsen theory. 

A more useful assumption is that $g=f$, so that we have the self coincidence situation. In this case,
if our path $c(t)$ does not close at $t=1$, we may reparametrize the cyclic homotopy so that $f_t=f$ for a
small interval ending at $t=1$, and let the path $c(t)$ connect up to the original starting point. Since we
have self coincidence everywhere in that small interval of $t$, we have $f(c(t)) = f(c(t))$ and we have obtained the desired loop.

Now  the coincidence number 
$$\Lambda_{f.g} = \roman{ind}(f, g, X)\tag 2.3$$
 and is given by the Lefschetz formula:
$$\Lambda_{f.g}=\Sigma_i(-1)^itr_i(f^*g^!) . \tag 2.4$$
Here, $f^*$ is the induced map of $f$ on integral cohomology $H^*(Y) \rightarrow H^*(X)$ and $g^!$ is the umkehr map (or Poincare duality map) from $H^*(X) \rightarrow H^*(Y)$ given by 
$D^{-1}_Yg_*D_X$ where $D_X: H^*(X)\to  H_*(X)$ is the Poincare duality isomorphism. Then
$tr_i(f^*g^!) $ denotes the trace of the composition on $H^i(X)$. An excellent description of Lefschetz coincidence theory is found in Chapter VI, section 14 of Glen Bredon's textbook [Bredon (1993)].

\medskip
\centerline{\bf 3. PROOFS OF THEOREMS}

\bigskip\bigskip

\demo {Proof of Theorem 1.4} 
In this case we have one map $f:M\to N$ and so the coincidence index 
$\roman{ind}(f, f, M)=\Lambda_{f,f}$ by (2.3). Now if we assume that $\roman{ind}(f, f, M)\neq 0$, we may apply Lemma 2.1 Now since we have a self coincidence, we may find the closed curve in the hypothesis of Lemma 2.2, and thus equation (2.2) holds and noting that $ \omega^f_*(\hat G)$ is trivial  we obtain
$$ \omega^f_*(\hat F) + f_*(\sigma_x) =  f_*(\sigma_x) \in \pi_1(N) .\tag 3.1$$
 Multiplying equation (3.1) on the right by the inverse of $f_*(\sigma_x)$ and gives us that 
 $$\omega^f_*(\hat F)= 1, \tag 3.2 $$
 that is it is trivial. Then the contrapositive gives us that a non-trivial 
$ \omega^f_*(\hat F)$  implies that $\Lambda_{f,f}=0$ \qed

\enddemo

Now Theorem 1.1 follows from the following lemma which computes the self intersection number.

\proclaim{Lemma 3.1} Let $f: M^n \rightarrow N^n$ be a map of closed oriented connected smooth manifolds of the same dimension. Then 
$$\Lambda_{f,f} = \text{deg}(f) \chi (N) \tag 3.3$$
\endproclaim

\demo{Proof} This result is actually mentioned in Bredon's textbook. We will give a short proof here for completeness. The following equation relates the umkehr map and cup products.
$$ f^!(f^*(a) \cup b)=a\cup f^!(b) \tag 3.4$$
Now letting $b=1\in H^0(M)$ we have $f^!(1)=\text{deg(f)}1\in H^0(N)$ and we obtain from (3.4)
$$ f^!(f^*(a))=\text{deg(f)}a \tag 3.5$$
Hence $tr_i(f^!f^*) = \text{deg}(f)\text{rank} (H^i(N))$. Of course, the alternating sums of the ranks gives the Euler Characteristic of $N$. So equation 2.4 combined with 
$tr_i(f^!f^*) = \text{deg}(f)\text{rank} (H^i(N))$ gives equation (3.3).\qed
\enddemo

\demo{proof of Theorem 1.1} Since $w^f_*$ is not trivial, Theorem 1.4 asserts that $\Lambda_{f,f}=0$
Hence $\Lambda_{f,f} = \text{deg}(f) \chi (N)=0$ by Lemma 1.3, so one of $ \text{deg}(f)$ or $ \chi (N)$ must equal zero. \qed
\enddemo

\demo{proof of Corollary 1.2} The space of maps of a CW complex  $X$ into an aspherical CW complex $K$ is aspherical itself, and its fundamental group is isomorphic to the centralizer of the of the image of
$f_*:\pi_1(X)\to \pi_1(K)$ in $ \pi_1(K)$. The evaluation map $\omega^f_*:\pi_1( \text{Map(X,K;f))}\to 
\pi_1(X)$ is injective. Special cases of these results appear in section III of [Gottlieb (1965)] and the full results are proved in ([Gottlieb (1969); see Lemma 2) and also in [Gottlieb (2004); lemma 2]. In the case at hand, $N$ is aspherical and so
if $\text{deg}(f)\times \chi(N) \neq 0$, then $\omega^f_*$ must be trivial which implies that the aspherical
space $\text{Map}(M,N;f)$ is simply connected, and hence since all the homotopy groups are trivial, so
$\text{Map}(M,N;f)$ is contractible. \qed
\enddemo

\demo{proof of Corollary 1.3}Let $N$ be a closed oriented connected aspherical manifold with 
Euler characteristic $ \chi(N) \neq 0$. Then every subgroup $G$ of finite index $k$ in $\pi_1(N)$ has a corresponding finite covering space $\tilde N$. The covering projection $p: \tilde N \to N$ is a map from a closed oriented smooth manifold of the same dimension as $N$ since $\tilde N$ is a finite covering,
and the degree of $p$ is $k \neq 0$. Hence if  $\chi (N) \neq 0$, the evaluation map 
$$\omega_*^p : \pi_1(\text{Map}(\tilde N,N;p) \rightarrow \pi_1(N)$$ must be trivial by Theorem 1.1. But 
$\omega_*^p$ is injective and maps $\pi_1(\text{Map}(\tilde N,N;p)$ onto the centralizer of $G$. Hence
the centralizer of $G$ in $\pi_1(N)$ must be trivial. \qed
\enddemo

\bigskip\bigskip

\centerline{\bf 4. FINAL THOUGHTS AND A CONJECTURE}

\bigskip\bigskip
Consider Corollary 1.2. If $f$ is the identity map, then $\text{deg}(f)=1$  and the hypothesis reduces to
$\chi(N)\neq 0$. Also the centralizer of the identity map $f$ is the center of $\pi_1(N)$ and we obtain the result: 
\proclaim{Corollary 4.1} If $N$ is an aspherical manifold so that $\chi(N)\neq 0$, then the center of  
$\pi_1(N)$ is trivial.
\endproclaim
This result was proved in [Gottlieb (1965); see Theorem IV.1] for $N$ a finite aspherical CW complex. The proof used Neilsen fixed point theory. Another proof with the same hypothesis but similar to the approach taken in Corollary 1.2, is in [Gottlieb (1990); see section 9]. This proof is very similar to the proof of Corollary 1.2 given here, where we used the index of a vector field to show that the appropriate
curve of zeros existed. 

Corollary 4.1 has had a remarkable history. The result was given an algebraic proof in [Stallings (1965)].
From that point on, the result become known as Gottlieb's theorem.

Gottlieb's theorem underwent a great generalization in [Rosset (1984)], in which the center of the fundamental group is replaced by a normal abelian subgroup. Later in [Cheeger--Gromov (1986)], the normal abelian subgroup was replaced by an infinite amenable subgroup. More  recently a sort of converse to Gottlieb's theorem was proved in [Hillman (1997)]. 

In addition to these strengthenings of the conclusion, the hypothesis of finite CW complex was somewhat relaxed. We know, however, that compactness cannot be totally eliminated below the homology level since
[Baumslag; Dyer; Heller (1980)] showed that there exists groups with acyclic homology with nontrivial
centers.

In view of the dramatic generalization of Gottlieb's theorem surveyed above, it is reasonable to ask if
Corollary 1.2 also is susceptible to generalization. The Euler Poincare number, of course, is a homotopy
invariant  for the homotopy category. We must generalize the degree of a map for continuous maps in general, instead of just for maps between oriented manifolds of the same dimension. We only need the absolute values of the degree. There are two generalizations which might work. They are found in [Gottlieb (1986)]:
Given a homomorphism $h: A \rightarrow B$ where $A$ and $B$ are abelian groups, a {\it transfer} for $h$ is a homopmorphism $\tau : B \rightarrow A$ so that $h \circ \tau = N$ where $N$ denotes multiplication on $B$ by the integer $N$. The set of integers associated to transfers for $h$ forms an ideal in the integers, and its positive generator is called the degree of $h$, $\text{deg}(h)$. If there is no positive generator, then $\text{deg}(h):=0$

Then we define the {\it degree} of a continuous map $f:X \to Y$ as the degree of the homomorphism 
   $f_*:H_*(X,\Bbb Z) \to H_*(Y, \Bbb Z)$.

\proclaim{Conjecture} Let $Y$ be an aspherical finite CW complex with $\chi (Y) \neq 0$, and let $\text{deg}(f) \neq 0$ where $f:X \to Y$. Then the centralizer of the image of $f_*$ in the fundamental group of
$Y$ is trivial. 
\endproclaim

The other approach is to define the codegree of $f$, denoted $\text{cdg}(f)$, by considering the induced
cohomology homomorphism $f^*: H^*(Y) \to H^*(X)$ and the associated transfers 
$\tau :  H^*(X)\to  H^*(Y)$ where $\tau \circ f^*= N$

\proclaim{Coconjecture} Let $Y$ be an aspherical finite CW complex with $\chi (Y) \neq 0$, and let $\text{cdg}(f) \neq 0$ where $f:X \to Y$. Then the centralizer of the image of $f_*$ in the fundamental group of
$Y$ is trivial. 
\endproclaim

Both degree and codegree generalize the absolute degree for maps between oriented manifolds of the same dimension. But degree and codegree do not agree in general. See [Gottlieb (1986)] for more details.

We can prove a special case of the conjecture for $\text{Map}(S^1,Y,f)$ with the condition on $Y$ that it is not a rationally acyclic. Note that 
$\text{Map}(S^1,Y)=: \Lambda (Y)$, the free loop space on $X$.

\proclaim{Proposition 4.2} For any map $f: S^1 \to Y$ for $Y$ a finite connected aspherical CW complex which has nontrivial rational homology above dimension zero, we have
$\text{deg}(f)\chi (Y) =0$.
\endproclaim

\demo{proof} Suppose that $\text{deg}(f)\neq 0$. Then there exists a transfer 
$\tau : H^*(S^1) \to H^*(Y)$ so that $f_* \circ \tau = N \neq 0$. This forces all the homology groups of dimension greater than 1 to be finite groups with each element having order dividing $N$. The 
homology groups are finitely generated since $Y$ is a finite complex. Now since some infinite homology group exists by hypothesis, in dimension greater than zero, it must be in dimension one. Also, for dimension 1, $H^1(Y)$ must have rank 1 by the transfer equation. Hence, rationally, $Y$ has the same homology as $S^1$, hence $\chi(Y)=0$. \qed
\enddemo

This result is a special case of the conjecture since for any map $f: S^1 \to Y$, the image of $f_*$ in the fundamental group $\pi_1(Y)$ has a nontrivial centralizer. Unfortunately, when $Y$ is a non simply connected finite aspherical 
acyclic finite complex, the conjectures are false, as pointed out by Thomas Schick. Schick also pointed out the necessity of the nonacyclic hypothesis on $Y$ in the above proposition.

So the above conjectures are false. but Schick and I would conjecture that there should be some reasonable conjecture. A more modest conjecture is:
\proclaim{Revised Conjecture} Let $X$ and $Y$ be Poincare duality complexes where  $Y$ is an aspherical Poincare duality complex with $\chi (Y) \neq 0$, and let $\text{deg}(f) \neq 0$ where $f:X \to Y$. Then the centralizer of the image of $f_*$ in the fundamental group of
$Y$ is trivial. 
\endproclaim

What follows is musings about the roles of various points of views which arose in this work:
There are three "regimes" in this subject. The mildest is the category of smooth closed manifolds and continuous maps. Here it appears that self coincidence Lefschetz numbers are subcases of Lefschetz
coincidence numbers; and fixed point Lefschetz numbers are special cases of coincidence numbers;
and fixed point numbers for maps homotopic to the identity are special cases of self coincidence numbers. But fixed point theory properly live in the more turbulent regime of the category of finite CW complexes and continuous maps. From the point of view of this new regime, fixed point theory does not generalize to coincidence theory.

Then by "accident", the results of fixed point theory are applied to aspherical spaces (ie. $K(\pi,1)$ and a result like Gottlieb's theorem is true. The category of $K(\pi,1))$'s and continuous maps is almost equivalent to the category of groups and homomorphisms, so Stallings' idea of finding a group theoretic
proof of Gottlieb's theorem led to steady weakening of the compactness condition and extensions of the subgroups from center to amenable leading eventually to the Cheeger--Gromov theorem. The Euler--Poincare number remains the same in this extension, as it does in so many other analogous situations. 

In the 1980's, I rediscovered a formula of [Morse (1929)], concerning the index of vector fields which made a big impression on me. It related the index of the vector field on a manifold with boundary to the index of a derived vector field on its boundary. It struck me that this formula itself provided an inductive definition on dimension of the index. I thought that all the properties of the index should be proved via this formula. That was the case, and I was able also to find quick proofs of famous theorems found in
Algebraic topology such as the Borsuk--Ulam theorem and the Brouwer fixed point theorem and 
Gottlieb's theorem as well as new results of the same geometric flavor. The intensity of my point of view can be judged from the series of papers this equation spawned: [Gottlieb (1990), (1988), (1986a), (1987)] and [Becker-Gottlieb (1991)]. 

Morse's formula depends only on continuity, the concept of pointing inside, the Euler--Poincare number and dimension. Out of this basic material flows integer invariants, the index of a particular zero, and a conservation law of zeroes moving under a homotopy of vector fields, which is reminiscent of particles in a cloud chamber. So I began to study physics.

I had a friendly debate with Albrecht Dold about the relationship of the index of a fixed point to the index of a zero of a vector field. In [Dold (1965)], the index of a fixed point composition of two maps is the same 
for either order of the composition. A corresponding result does not exist for indices of zeros of vector fields.
On the other hand, the index of a vector field $V$ is equal to the index of $-V$ for even dimensions and equal to its negative for odd dimensions.  An important fact which has no immediate correspondence for the fixed point index.

I gave my student the problem to find a Morse formula for fixed point indices, see [Benjamin (1990)].
After more than 10 years we published the formula in [Benjamin and  Gottlieb (2006)]. An attempt was made in [Benjamin (1990)] to find the index of a path field for topological manifolds. But she couldn't prove the fact that $\text{ind}(V)= \pm \text{ind}(V)$, except for the case of a smooth manifold, which gave nothing new. Path fields were invented in [Nash (1955] and developed further in [Fadell (1965)] and 
[Brown (1965)].

Despite my hopes, the Morse formula for fixed point indices was not as productive as the Morse formula for vector fields.

The point of the last paragraphs was to show how problematic the notion of the correct point of view is.

\vfill\eject

\eject
\enddocument
\bye